# AREAS OF SPHERICAL AND HYPERBOLIC TRIANGLES IN TERMS OF THEIR MIDPOINTS

G.M. Tuynman

ABSTRACT. Let $M$ be either the 2-sphere $\mathbf{S}^2 \subset \mathbf{R}^3$ or the hyperbolic plane $\mathbf{H}^2 \subset \mathbf{R}^3$. If $\Delta(abc)$ is a geodesic triangle on $M$ with corners at $a, b, c \in M$, we denote by $\alpha, \beta, \gamma \in M$ the midpoints of their sides. If $\Omega$ denotes the oriented area of this triangle on $M$, it satisfies the relations :

$$(0) \qquad \sin(\Omega/2) = \frac{\det_3(abc)}{\sqrt{2(1 + \langle a, b \rangle)(1 + \langle b, c \rangle)(1 + \langle c, a \rangle)}} = \det_3(\alpha\beta\gamma) \ ,$$

where $\langle \ , \ \rangle$ denotes the Euclidean scalar product for $M = \mathbf{S}^2$ and the Lorentzian scalar product for $M = \mathbf{H}^2$. On the hyperbolic plane one should always take the solution with $|\Omega/2| < \pi/2$. On the sphere, singular cases excepted, a straightforward procedure tells us which solution of this equation is the correct one.

## General introduction

Three points $a$, $b$, and $c$ in the Euclidean plane $\mathbf{R}^2$ determine a unique (geodesic) triangle $\Delta(abc)$. Moreover, the three midpoints $\alpha$, $\beta$, and $\gamma$ of the three sides $\overline{bc}$, $\overline{ca}$, and $\overline{ab}$ respectively of $\Delta$ also determine this triangle (the side $\overline{ab}$ passes through $\gamma$ and is parallel to $\overline{\alpha\beta}$). If $\Omega$ denotes the oriented area of the triangle $\Delta$, it is given by the formulæ

$$\Omega = \tfrac{1}{2} \det_2(b - a, c - a) = 2 \det_2(\beta - \alpha, \gamma - \alpha) \ ,$$

where $\det_2$ denotes the $2 \times 2$ determinant. If we view the Euclidean plane $\mathbf{R}^2$ as the the subset of $\mathbf{R}^3$ with the third coordinate equal to 1, these formulæ can be written in the more symmetric form

$$\Omega/2 = \tfrac{1}{4} \det_3 \begin{pmatrix} a & b & c \\ 1 & 1 & 1 \end{pmatrix} = \det_3 \begin{pmatrix} \alpha & \beta & \gamma \\ 1 & 1 & 1 \end{pmatrix} \ .$$

This simple situation changes when we replace the Euclidean plane $\mathbf{R}^2$ by either the 2-sphere $\mathbf{S}^2 \subset \mathbf{R}^3$ or the hyperbolic plane $\mathbf{H}^2$ seen as a sheet of the 2-sheeted hyperboloid in $\mathbf{R}^3$. The purpose of this paper is to find the equivalent formulæ for the sperical and hyperbolic cases.

The 2-sphere

**Introduction**

If $a$ and $b$ are two distinct points on $\mathbf{S}^2$, they determine a unique great circle (geodesic), provided they are not antipodal. But since this geodesic is closed, we have a choice for the geodesic segment $\overline{ab}$. As a consequence, if $c$ is a third point on $\mathbf{S}^2$, there is more than one triangle with $a$, $b$, and $c$ as corners (see figure 1, in which $a$ is the north pole and $b$ and $c$ lie on the equator).

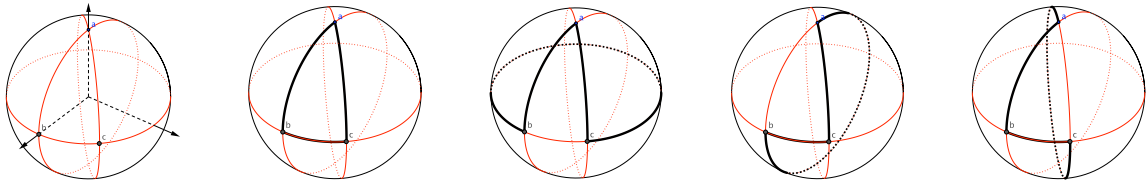

At a first glance one might think that there should be 8 possibilities, since for each side we have 2 choices. However, for each side one of the choices has a length less than $\pi$ and the other has length greater than $\pi$ (a great circle has length $2\pi$). Now if two sides have both length greater than $\pi$, these sides must intersect in antipodal points. Since this prohibits the construction of a well defined triangle, only one side can have length greater than $\pi$. It follows that there are "only" 4 possible choices for the triangle $\Delta(abc)$. Looking at the area of these various possibilities, it is obvious that the choice with three sides shorter than $\pi$ together with a choice for which one side is greater than $\pi$ fill out a hemisphere (see figure 1). It follows that the area $\Omega$ of the triangle with all sides less than $\pi$ and the area $\Omega'$ of the triangle with one side greater than $\pi$ satisfy the relation $\Omega + \Omega' = 2\pi =$ the area of a hemisphere (note that this formula is compatible with formula (0)). On the other hand, specifying the midpoint of a side fixes the choice described above. This leads to the idea that specifying the three midpoints $\alpha$, $\beta$, and $\gamma$ (defined as in the planar case) completely determines the triangle $\Delta$ on $\mathbf{S}^2$, and thus its area. Our purpose is to prove that this is indeed the case, except in some well defined singular cases.

One final word on the definition of area on the sphere. Since the sphere is orientable, we can speak of the oriented area (after having chosen an orientation). Now if we have three geodesic segments on the sphere forming a triangle, we still have to specify its "interior" in order to compute the area. If we have an oriented triangle, it is immediate that the *difference* of the *oriented* areas of the two choices for the interior equals the area of the sphere, i.e., $4\pi$. It follows that if we only look at the oriented area *modulo* $4\pi$, then the choice of interior is irrelevant.

**The analysis**

In this section we will denote a point in $\mathbf{R}^3 = \mathbf{C} \times \mathbf{R}$ by the coordinates $(w, h)$, with $w \in \mathbf{C}$ and $h \in \mathbf{R}$. A point $(w, h)$ lies on the unit sphere $\mathbf{S}^2$ if and only if $|w|^2 + h^2 = 1$. On $\mathbf{S}^2$ we will also use polar coordinates $(\theta, \varphi)$, or stereographic projection (from the north pole) coordinates $z \in \mathbf{C}$. Between these coordinate systems on $\mathbf{S}^2$ we have the following relations :

$$z = \operatorname{cotg}(\frac{\theta}{2})\, e^{i\varphi} = \frac{w}{1-h} \quad,\quad w = \sin\theta\, e^{i\varphi} = \frac{2z}{|z|^2+1} \quad,\quad h = \cos\theta = \frac{|z|^2-1}{|z|^2+1}\ .$$

Let $a$, $b$, $c$ be three points on the unit sphere $\mathbf{S}^2 \subset \mathbf{R}^3$. If there are no antipodal points among these three points, then there exists a unique triangle $\Delta$ with these points as corners for which all sides have length less than $\pi$. Our first goal will be to determine the area $\Omega$ of this triangle. Since the oriented area is invariant under rotations (the group $SO(3)$), we may assume without loss of generality that $a$ is the north pole. In spherical coordinates we thus have $w_b = \cos\theta_b\, e^{i\varphi_b}$, $w_a = 0$, and $w_c = \cos\theta_c\, e^{i\varphi_c}$ with $\theta_b, \theta_c \in [0, \pi)$; we may assume without loss of generality that $\Phi = \varphi_c - \varphi_b$ belongs to $(-\pi, \pi)$. The geodesic $\overline{ab}$ is described by $\varphi = \varphi_b$, $\theta \in (0, \theta_b)$ and the geodesic $\overline{ca}$ by $\varphi = \varphi_c$, $\theta \in (0, \theta_c)$. A point $g$ on the geodesic $\overline{bc}$ is described by the condition that $a$, $c$, and $g$ lie in the same plane, i.e., $\det(acg) = 0$. This gives the relation

$$\cotg \theta_g = A \sin\varphi_g + B \cos\varphi_g = \frac{\sin(\varphi_c - \varphi_g)}{\sin \Phi} \cotg \theta_b + \frac{\sin(\varphi_g - \varphi_b)}{\sin \Phi} \cotg \theta_c$$

with

$$A = \frac{\cotg \theta_c \cos\varphi_b - \cotg \theta_b \cos\varphi_c}{\sin \Phi} \quad , \quad B = \frac{\cotg \theta_b \sin\varphi_c - \cotg \theta_c \sin\varphi_b}{\sin \Phi} .$$

The oriented area element on $\mathbf{S}^2$ is given by $\sin\theta\, d\theta \wedge d\varphi$, and thus the oriented area $\Omega$ of the triangle $\Delta(abc)$ is given by

$$\Omega = \int_{\varphi_b}^{\varphi_c} \left( \int_0^{\theta_g(\varphi)} \sin\theta\, d\theta \right) d\varphi = \Phi - C_{bc} - C_{cb} ,$$

$$C_{jk} = \arcsin \left( \frac{|\sin \Phi|}{\sin \Phi} \cdot \frac{\cotg \theta_j - \cotg \theta_k \cos \Phi}{\sqrt{\cotg^2 \theta_j + \cotg^2 \theta_k - 2\cotg \theta_j\, \cotg \theta_k \cos \Phi + \sin^2 \Phi}} \right)$$

We then introduce $\alpha_{jk} = (\cotg \theta_j - \cotg \theta_k \cos \Phi) - i \sin \Phi / \sin \theta_k$ to obtain

$$-i\, e^{iC_{jk}} = \sin C_{jk} - i \cos C_{jk} = \frac{|\sin \Phi|}{\sin \Phi} \cdot \frac{\alpha_{jk}}{|\alpha_{jk}|} ,$$

where we used that $\cos C_{jk} = +\sqrt{1 - \sin^2 C_{jk}}$ because the arcsine takes its values in $(-\pi/2, \pi/2)$. Another elementary computation gives

$$\alpha_{cb} = \frac{e^{i\varphi_b}(\bar{z}_b z_c + 1)(\bar{z}_c - \bar{z}_b)}{2|z_b z_c|} \quad \& \quad \alpha_{bc} = \frac{e^{-i\varphi_c}(\bar{z}_b z_c + 1)(z_b - z_c)}{2|z_b z_c|} .$$

Since $e^{i\varphi} = z/|z|$, one deduces from the above that

$$e^{i\Omega} = e^{i(\Phi - C_{bc} - C_{cb})} = e^{i\Phi} \cdot \overline{e^{iC_{bc}}} \cdot \overline{e^{iC_{cb}}} = \left( \frac{\bar{z}_b z_c(z_b \bar{z}_c + 1)}{|\bar{z}_b z_c(z_b \bar{z}_c + 1)|} \right)^2 .$$

From this formula one deduces immediately that $\Omega \mod 4\pi$ and $2\operatorname{Arg}(\bar{z}_b z_c(z_b \bar{z}_c + 1))$ coincide modulo $2\pi$. The argument of Arg is zero if and only if $z_b = 0$ or $z_c = 0$, which means $b$ or $c$ the south pole, or $z_b \bar{z}_c + 1 = 0$, which means that $b$ and $c$ are antipodal. Since we have excluded these possibilities from our discussion, we conclude that $\Omega \mod 4\pi$ and $2\operatorname{Arg}(\bar{z}_b z_c(z_b \bar{z}_c + 1))$ are both continuous functions of $z_b$ and $z_c$ with values in $\mathbf{R}/4\pi\mathbf{R}$. Hence they must coincide. In order to express this result in a way which does not presuppose that $a$ is the north pole, we note that we have the following identities when *a is* the north pole:

$$
(1) \quad \langle a, b\rangle_E = h_b \;,\quad \langle a, c\rangle_E = h_c \;,\quad \langle b, c\rangle_E = h_b h_c + \Re(\overline{w}_b w_c)
$$
$$
\det{}_3(abc) = \Im(\overline{w}_b w_c) \;,
$$

where $\langle \cdot, \cdot \rangle_E$ denotes the usual Euclidean scalar product in $\mathbf{R}^3$ and $\det_3$ the $3 \times 3$ determinant. Since $\bar{z}_b z_c(z_b \bar{z}_c + 1) = \dfrac{1 + h_b + h_c + h_b h_c + \overline{w}_b w_c}{(1-h_b)(1-h_c)}$, we obtain the following lemma.

**Lemma.** *The oriented area modulo $4\pi$ of the spherical triangle $\Delta(abc)$ with all sides shorter than $\pi$ is given by*

$$
\begin{aligned}
\Omega &= 2\operatorname{Arg}\!\left(\bar{z}_b z_c(z_b \bar{z}_c + 1)\right) \qquad \text{if } a \text{ is the north pole} \\
&= 2\operatorname{Arg}\!\left(1 + \langle a, b\rangle_E + \langle b, c\rangle_E + \langle a, c\rangle_E + i\det{}_3(abc)\right) \;.
\end{aligned}
$$

We now turn our attention back to the area of a triangle determined by the midpoints of its sides. Without loss of generality we may assume that $\overline{ab}$ is the longest side and that its midpoint $\gamma$ is the north pole, which implies in particular that $z_a = -z_b$, $w_a = -w_b$, and $h_a = h_b$. We thus can cut the triangle $\Delta(abc)$ into two triangles $\Delta(\gamma bc)$ and $\Delta(\gamma ca)$ which both satisfy the condition that all sides are shorter than $\pi$ and that one corner ($\gamma$) is the north pole. Applying the previous lemma thus gives us for the surface $\Omega$ of the triangle $\Delta(abc)$ the formula:

$$
\begin{aligned}
\Omega &= 2\operatorname{Arg}\!\left(\bar{z}_b z_c(z_b \bar{z}_c + 1)\right) + 2\operatorname{Arg}\!\left(\bar{z}_c(-z_b)(z_c(-\bar{z}_b) + 1)\right) \\
(2) \qquad &= 2\operatorname{Arg}\!\left((\bar{z}_b z_c - 1)(z_b \bar{z}_c + 1)\right) \;,
\end{aligned}
$$

where the last equality is justified by the same type of continuity argument as the one which led to the above lemma. This formula gives the area of the triangle in terms of two of its corners ($b$ and $c$) and one midpoint (the north pole $\gamma$). Our next step will be to compute the remaining midpoints in terms of the corners, and then to invert this map. That way we will get the area in terms of the three midpoints.

Since the lengths of the two segments $\overline{bc}$ and $\overline{ca}$ are less than $\pi$, their midpoints $\alpha$ and $\beta$ are the normalized averages of the endpoints. Thus, if we introduce $\sigma_\alpha = w_c + w_b$, $\sigma_\beta = w_c - w_b$, and $\tau = h_c + h_b$, then $w_\epsilon = \sigma_\epsilon/\sqrt{|\sigma_\epsilon|^2 + \tau^2}$ and $h_\epsilon = \tau/\sqrt{|\sigma_\epsilon|^2 + \tau^2}$, $\epsilon = \alpha, \beta$ (remember that $w_a = -w_b$ and $h_a = h_b$). Using formula (1) and the relation $|w|^2 + h^2 = 1$, we find $\langle \alpha, \beta\rangle_E = 2\tau h_b((|\sigma_\alpha|^2 + \tau^2)(|\sigma_\beta|^2 + \tau^2))^{-1/2}$. This proves the first part of the following lemma; the corollary follows directly from formula (1).

**Lemma.** *We always have $\operatorname{sign} h_\alpha = \operatorname{sign} h_\beta = \operatorname{sign} \tau$; in case $h_\alpha = h_\beta = 0$, we have $\langle \alpha, \beta\rangle_E = 0$.*

**Corollary.** *If $\overline{ab}$ is the longest side of the triangle $\Delta(abc)$, the inner products $\langle \alpha, \gamma \rangle_E$ and $\langle \beta, \gamma \rangle_E$ have the same sign; if they are zero the vectors $\alpha$, $\beta$, and $\gamma$ are orthogonal. Moreover, $\langle \alpha, \beta \rangle_E$ has the same sign as the other two inner products if and only if the length of $\overline{ab}$ is less than $\pi$.*

The next task is to invert the application $(z_b, z_c) \mapsto (z_\alpha, z_\beta)$, for which we obviously have to take the constraints of the previous lemma into account, i.e., $h_\alpha$ and $h_\beta$ must have the same sign. From its definition we find $w_\epsilon/h_\epsilon = \sigma_\epsilon/\tau$, which leads us to introduce $s_c = \frac{1}{2}(w_\beta/h_\beta + w_\alpha/h_\alpha) = w_c/\tau$ and $s_b = \frac{1}{2}(w_\alpha/h_\alpha - w_\beta/h_\beta) = w_b/\tau$. We thus have $w_j = (h_b + h_c)s_j$ and $|w_j|^2 + h_j^2 = 1$ ($j = b, c$). Solving these equations for $h_b$ and $h_c$ involves an equation of degree 2, i.e., a sign $\eta = \pm$, and one finds (using $h_b + h_c = \tau \neq 0$) :

$$h_b = \frac{\eta(1 + |s_c|^2 - |s_b|^2)}{\sqrt{(1 + |s_c|^2 + |s_b|^2)^2 - 4|s_b s_c|^2}} \quad , \quad \tau = \frac{2\eta}{\sqrt{(1 + |s_c|^2 + |s_b|^2)^2 - 4|s_b s_c|^2}} \ .$$

Since $\eta = \operatorname{sign} \tau = \operatorname{sign} h_\alpha = \operatorname{sign} h_\beta$, this sign is well determined and thus we find a unique solution for $(w_j, h_j)$, $j = a, b$ in terms of $(w_\epsilon, h_\epsilon)$, $\epsilon = \alpha, \beta$ :

$$h_b = \eta \frac{h_\alpha h_\beta + \Re(w_\beta \overline{w}_\alpha)}{\sqrt{1 - (\Im(w_\beta \overline{w}_\alpha))^2}} \quad , \quad h_c = \eta \frac{h_\alpha h_\beta - \Re(w_\beta \overline{w}_\alpha)}{\sqrt{1 - (\Im(w_\beta \overline{w}_\alpha))^2}}$$

$$w_b = \eta \frac{w_\alpha h_\beta - w_\beta h_\alpha}{\sqrt{1 - (\Im(w_\beta \overline{w}_\alpha))^2}} \quad , \quad w_c = \eta \frac{w_\beta h_\alpha + w_\alpha h_\beta}{\sqrt{1 - (\Im(w_\beta \overline{w}_\alpha))^2}} \ .$$

We see that there exists a unique inverse provided $(\Im(w_\beta \overline{w}_\alpha))^2$ is not 1, which occurs if and only if $h_\beta = h_\alpha = 0$ (remember the relation $|w|^2 + h^2 = 1$). Converting this into $z = w/(1-h)$ coordinates yields :

$$(\overline{z}_b z_c - 1)(z_b \overline{z}_c + 1) = \frac{4 h_\alpha h_\beta \left( \eta \sqrt{1 - (\Im(w_\beta \overline{w}_\alpha))^2} + i\Im(w_\beta \overline{w}_\alpha) \right)}{(1 - h_b)(1 - h_c)(1 - (\Im(w_\beta \overline{w}_\alpha))^2)} \ ,$$

and thus we find from formula (2) for the area :

$$\Omega = 2 \operatorname{Arg} \left( \eta \sqrt{1 - (\Im(w_\beta \overline{w}_\alpha))^2} + i \Im(w_\beta \overline{w}_\alpha) \right) \ .$$

Since the third midpoint is the north pole, we can apply formula (1) to show that $\Im(w_\beta \overline{w}_\alpha) = \det_3(\alpha\beta\gamma)$. Using the corollary, we thus have proven the following proposition.

**Proposition.** *Let $\alpha$, $\beta$, and $\gamma$ be three points on the unit sphere such that none of the inner products is zero. Then the oriented area $\Omega$ of the triangle $\Delta(abc)$ whose midpoints of the three sides are the given points $\alpha$, $\beta$, and $\gamma$ is given by*

$$e^{i\Omega/2} = \eta \sqrt{1 - \det(\alpha\beta\gamma)^2} + i \det(\alpha\beta\gamma) \ ,$$

*where $\eta$ is the sign of the majority of the signs of the three inner products $\langle \alpha, \beta \rangle_E$, $\langle \beta, \gamma \rangle_E$, and $\langle \gamma, \alpha \rangle_E$.*

The obvious next question is what happens if one or more of the inner products is zero. In the next section we will show that, if one or two inner products are zero, then necessarily there are antipodal points among the corners $a$, $b$, and $c$. The next figures then show that these cases do not allow us to determine the triangle $\Delta(abc)$ from its midpoints $\alpha$, $\beta$, and $\gamma$ (one side is the first part of another side of length greater than $\pi$). Moreover, these figures show that it is not possible to make a consistent choice for the area $\Omega$ such that it becomes a continuous function of the midpoints.

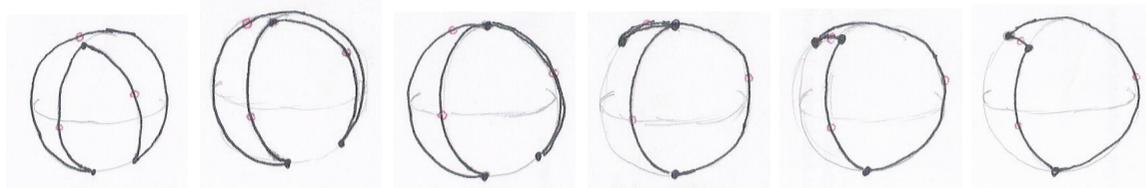

Finally, all three inner products are zero if only if the three points $\alpha$, $\beta$, and $\gamma$ form an orthonormal frame, i.e., when $\det_3(\alpha\beta\gamma)^2 = 1$. It follows that our proposition tells us that the area satisfies $e^{i\Omega/2} = i\det(\alpha\beta\gamma) = \pm i$, i.e., the area is well defined. On the other hand, we will see in the next section that in this case the triangle $\Delta(abc)$ is completely undetermined. A geometric way to see this happening is as follows. Suppose there exist two triangles $\Delta(abc)$ and $\Delta(a'b'c')$ that have the same midpoints of their sides. Given that $\gamma$ is the midpoint of both $\overline{ab}$ and $\overline{a'b'}$, it follows that the move $a \to a'$ completely determines the move $b \to b'$. In turn, this move determines the move $c \to c'$, which in turn redetermines the move $a \to a'$. But this redefined move is turned with respect to the initial move over the angle $\angle abc + \angle bca + \angle cab$. It follows that the sum of the three angles in the triangle $\Delta(abc)$ must be a multiple of $2\pi$. On the other hand, the excess angle, i.e., the sum of the three angles minus $\pi$, is the area $\Omega$. We thus find that if there exist different (in a generic sense) triangles with the same midpoints, then necessarily $\Omega = \pm\pi$ modulo $4\pi$, i.e., $e^{i\Omega/2} = \pm i$.

**A qualitative analysis**

We will denote by $\mathcal{R}_v^\phi$ a rotation in $\mathbf{R}^3$ over an angle $\phi$ around the axis $v \in \mathbf{R}^3$. For any $\lambda > 0$ we have $\mathcal{R}_{\lambda v}^\phi = \mathcal{R}_v^\phi$, but for $\lambda < 0$ we have $\mathcal{R}_{\lambda v}^\phi = \mathcal{R}_v^{\phi+\pi}$. For future use we recall that $\mathcal{R}_v^\pi \circ \mathcal{R}_w^\pi = \mathcal{R}_{v \wedge w}^\phi$, where $v \wedge w$ denotes the vector product of the two vectors $v, w \in \mathbf{R}^3$, and where $\phi$ is twice the angle between $v$ and $w$.

Since $\alpha$ is the midpoint of $\overline{bc}$, we have $\mathcal{R}_\alpha^\pi(b) = c$. Similarly we have $\mathcal{R}_\beta^\pi(c) = a$ and $\mathcal{R}_\gamma^\pi(a) = b$. We conclude that $b$ is a fixed point of the rotation $\mathcal{R}_\gamma^\pi \circ \mathcal{R}_\beta^\pi \circ \mathcal{R}_\alpha^\pi$. Now, unless it is the identity, a rotation has a unique axis of rotation, and thus only two (antipodal) fixed points on the unit sphere $\mathbf{S}^2$. If $b$ denotes one of these two fixed points, the points $a$ and $c$ are uniquely determined by $c = \mathcal{R}_\alpha^\pi(b)$ and $a = \mathcal{R}_\gamma^\pi(b) = (\mathcal{R}_\beta^\pi \circ \mathcal{R}_\alpha^\pi)(b)$. The segment $\overline{bc}$ is uniquely determined by the fact that $\alpha$ is its midpoint, even if $b$ and $c$ are opposite points on the sphere. Still under the assumption that $\mathcal{R}_\gamma^\pi \circ \mathcal{R}_\beta^\pi \circ \mathcal{R}_\alpha^\pi$ is not the identity, we thus conclude that there exist exactly two "triangles" $\Delta(abc)$ determined by the midpoints $\alpha$, $\beta$, and $\gamma$. However, changing from one solution to the other, we will change the lengths $\ell$ of the sides to their complement $2\pi - \ell$. Since only one side can be longer than $\pi$, we conclude that the midpoints $\alpha$, $\beta$, and $\gamma$ determine a unique triangle $\Delta(abc)$, provided that there are no antipodal points among the solutions for $a$, $b$, and $c$ (which would give $\ell = 2\pi - \ell = \pi$), and provided that $\mathcal{R}_\gamma^\pi \circ \mathcal{R}_\beta^\pi \circ \mathcal{R}_\alpha^\pi$ is not the identity.

**Lemma.** *The rotation $\mathcal{R}^\pi_\gamma \circ \mathcal{R}^\pi_\beta \circ \mathcal{R}^\pi_\alpha$ is the identity if and only if $\alpha$, $\beta$, and $\gamma$ form an orthonormal basis of $\mathbf{R}^3$, i.e., $\det_3(\alpha\beta\gamma)^2 = 1$.*

*Proof.* The rotation $\mathcal{R}^\pi_\gamma \circ \mathcal{R}^\pi_\beta \circ \mathcal{R}^\pi_\alpha$ is the identity if and only if $\mathcal{R}^\pi_\gamma = \mathcal{R}^\phi_{\alpha\wedge\beta}$ with $\phi$ twice the angle between $\alpha$ and $\beta$. This is the case if and only if $\gamma$ is parallel to $\alpha \wedge \beta$ and $\phi = \pi$, which happens if and only if the angle between $\alpha$ and $\beta$ is $\pi/2$ and $\gamma$ orthogonal to the plane spanned by $\alpha$ and $\beta$. $\boxed{QED}$

If $\det_3(\alpha\beta\gamma)^2 = 1$, then $\mathcal{R}^\pi_\gamma \circ \mathcal{R}^\pi_\beta \circ \mathcal{R}^\pi_\alpha$ is the identity and thus we can choose the fixed point $b$ freely (modulo the condition that the three points $a$, $b$, and $c$ form a genuine triangle). At the end of the previous section we have argued that in that case the area must be $\pm\pi$ depending upon the orientation. This is compatible with our formula $\sin(\Omega/2) = \det_3(\alpha\beta\gamma)$.

Having determined when $\mathcal{R}^\pi_\gamma \circ \mathcal{R}^\pi_\beta \circ \mathcal{R}^\pi_\alpha$ is the identity, we now turn our attention to the case where there are antipodal points among the solutions for $a$, $b$, and $c$. If $b$ and $c$ are antipodal, then certainly the lines $\mathbf{R}b$ and $\mathbf{R}c$ coincide. But these are the axes of the rotations $\mathcal{R}^\pi_\gamma \circ \mathcal{R}^\pi_\beta \circ \mathcal{R}^\pi_\alpha$ (of which $b$ is a fixed point) and $\mathcal{R}^\pi_\alpha \circ \mathcal{R}^\pi_\gamma \circ \mathcal{R}^\pi_\beta$ (of which $c$ is a fixed point).

**Lemma.** *If $\mathcal{R}^\pi_\gamma \circ \mathcal{R}^\pi_\beta \circ \mathcal{R}^\pi_\alpha$ is not the identity, then $A = \mathcal{R}^\pi_\gamma \circ \mathcal{R}^\pi_\beta \circ \mathcal{R}^\pi_\alpha$ and $B = \mathcal{R}^\pi_\alpha \circ \mathcal{R}^\pi_\gamma \circ \mathcal{R}^\pi_\beta$ have the same axis of rotation if and only if either*

(i) $\mathbf{R}\beta = \mathbf{R}\gamma$, *in which case $b = c = \pm\alpha$;*
(ii) $\langle \gamma, \beta \rangle_E = 0$, $\det(\alpha\beta\gamma)^2 \neq 1$, *in which case $c = -b$;*
(iii) $\langle \alpha, \beta \rangle_E = 0 = \langle \alpha, \gamma \rangle_E$, $\langle \gamma, \beta \rangle_E \neq 0$, *in which case $-a = b = c = \pm\alpha$.*

*Proof.* $\mathbf{R}\beta = \mathbf{R}\gamma$ if and only if $\beta$ and $\gamma$ are parallel if and only if $\beta \wedge \gamma = 0$, and then $\mathcal{R}^\pi_\gamma \circ \mathcal{R}^\pi_\beta$ is the identity. It follows that $A = B = \mathcal{R}^\pi_\alpha$ and thus $\pm\alpha = b = \mathcal{R}^\pi_\alpha(b) = c$.

So let us suppose $\beta\wedge\gamma \neq 0$, which implies that $\mathcal{R}^\pi_\gamma \circ \mathcal{R}^\pi_\beta = \mathcal{R}^\phi_{\gamma\wedge\beta}$ is not the identity ($\phi$, twice the angle between $\beta$ and $\gamma$, is different from $0$ and $2\pi$). Since $A$ and $B$ are rotations different from the identity, having the same axis means that there exists a non-zero vector $b$ such that $Ab = b$ and $Bb = b$. Since $\mathcal{R}^\pi_\alpha \circ A = B \circ \mathcal{R}^\pi_\alpha$, this implies that $B(\mathcal{R}^\pi_\alpha(b)) = \mathcal{R}^\pi_\alpha(b)$, and thus that $\mathcal{R}^\pi_\alpha(b) = \lambda b$ with $\lambda = \pm 1$ (because a rotation can only have $\pm 1$ as eigenvalues). We then compute : $b = Ab = (\mathcal{R}^\phi_{\gamma\wedge\beta})(\mathcal{R}^\pi_\alpha(b)) = \lambda(\mathcal{R}^\phi_{\gamma\wedge\beta})(b)$, and thus $\mathcal{R}^\phi_{\gamma\wedge\beta}(b) = \lambda b$ as well as $\mathcal{R}^\pi_\alpha(b) = \lambda b$.

If $\lambda = 1$, $\mathbf{R}b$ must be the axis of rotation, and thus in particular $\alpha$ and $\gamma \wedge \beta$ must be parallel, i.e., $\langle \alpha, \beta \rangle_E = 0$ and $\langle \alpha, \gamma \rangle_E = 0$; in this situation $\langle \gamma, \beta \rangle_E = 0$ is excluded by the condition that $A$ is not the identity. $\mathbf{R}b$ being the axis of rotation of $\mathcal{R}^\pi_\alpha$, we have $\pm\alpha = b = \mathcal{R}^\pi_\alpha(b) = c$. Moreover, $\alpha$ being orthogonal to $\gamma$, it follows that $a = \mathcal{R}^\pi_\gamma(b) = -b$.

If $\lambda = -1$, $b$ must be perpendicular to the axis of rotation of $\mathcal{R}^\pi_\alpha$, and necessarily the rotation angle must be $\pi$. This implies on the one hand that $c = \mathcal{R}^\pi_\alpha(b) = -b$, and on the other hand that $\phi = \pi$, i.e., $\langle \gamma, \beta \rangle_E = 0$. Since $A = \mathcal{R}^\pi_{\gamma\wedge\beta} \circ \mathcal{R}^\pi_\alpha = \mathcal{R}^\psi_{(\gamma\wedge\beta)\wedge\alpha}$ is not the identity, $(\gamma \wedge \beta) \wedge \alpha \neq 0$, which is equivalent to $\det(\alpha\beta\gamma)^2 \neq 1$. $\boxed{QED}$

**Corollary.** *The only cases in which the midpoints $\alpha$, $\beta$, and $\gamma$ (possibly) determine antipodal points among the $a$, $b$, and $c$ are :*

(i) *if the three inner products $\langle \alpha, \beta \rangle_E$, $\langle \beta, \gamma \rangle_E$, and $\langle \alpha, \gamma \rangle_E$ are all zero, in which case the corners $a$, $b$, and $c$ are undetermined;*

(ii) *if two out of the three inner products are zero, in which case two corners coincide and the third is antipodal;*
(iii) *if one out of the three inner products is zero, in which case two corners are antipodal.*

## The hyperbolic plane

**Introduction**

We view the hyperbolic plane $\mathbf{H}^2$ as the upper sheet of the 2-sheeted hyperboloid in $\mathbf{R}^3$. As before we use coordinates $(w, h) \in \mathbf{C} \times \mathbf{R} = \mathbf{R}^3$ and thus $\mathbf{H}^2$ is determined by the equations $h > 0$ and $h^2 - |w|^2 = 1$. On $\mathbf{R}^3$ we also introduce the Lorentzian scalar product $\langle\,,\,\rangle_L$ defined as

$$\langle (w, h), (w', h') \rangle_L = hh' - \Re(\overline{w}w') \ .$$

This Lorentzian scalar product induces a metric and an oriented area element on $\mathbf{H}^2$, which in hyperbolic coordinates $(\theta, \varphi)$ (see below) is given as $\sinh\theta\, d\theta \wedge d\varphi$. In this realization of $\mathbf{H}^2$, the geodesics are just the intersections of $\mathbf{H}^2$ with planes in $\mathbf{R}^3$ passing through the origin. Contrary to the spherical case, two points on $\mathbf{H}^2$ determine a unique geodesic segment between them and thus a unique midpoint. Using that $SO(2,1)$ (the subgroup of $SL(3,\mathbf{R})$ preserving $\langle\,,\,\rangle_L$) acts transitively on $\mathbf{H}^2$, it is easy to show that the midpoint $\gamma$ on the geodesic between $a$ and $b$ in $\mathbf{H}^2$ is given as the normalized average, i.e., $\gamma = (a+b) \cdot (\langle a+b, a+b\rangle_L)^{-1/2}$.

**The analysis**

The computations to find the area of a hyperbolic triangle follow closely those done for the sphere. As said, on $\mathbf{R}^3$ we use coordinates $(w,h) \in \mathbf{C} \times \mathbf{R}$. We also use hyperbolic coordinates $(\theta, \varphi)$ with $\theta \in [0, \infty)$, or stereographic projection (from the south pole $(0+0i, -1)$) coordinates $z \in \mathbf{C}$, $|z|^2 < 1$. Between these coordinate systems on $\mathbf{H}^2$ we have the following relations

$$z = \tanh(\frac{\theta}{2})\, \mathrm{e}^{i\varphi} = \frac{w}{1+h} \quad , \quad w = \sinh\theta\, \mathrm{e}^{i\varphi} = \frac{2z}{1-|z|^2} \quad , \quad h = \cosh\theta = \frac{1+|z|^2}{1-|z|^2} \ .$$

Let $a$, $b$, $c$ be three points on the hyperbolic plane $\mathbf{H}^2 \subset \mathbf{R}^3$, and let $\Delta$ be the unique triangle with these points as corners. Our first goal will be to determine the area $\Omega$ of this triangle. Since the area is invariant under $SO(2,1)$, we may assume without loss of generality that $a$ is the north pole $(0+0i, 1)$. In hyperbolic coordinates we thus have $w_b = \cosh\theta_b\, \mathrm{e}^{i\varphi_b}$, $w_a = 0$, and $w_c = \cosh\theta_c\, \mathrm{e}^{i\varphi_c}$; we may assume without loss of generality that $\Phi = \varphi_c - \varphi_b$ belongs to $(-\pi, \pi)$. If $\Phi = \pm\pi$, $b$ and $c$ are opposite, and thus we have a degenerate triangle with area zero. The geodesic $\overline{ab}$ is described by $\varphi = \varphi_b$, $\theta \in (0, \theta_b)$ and the geodesic $\overline{ca}$ by $\varphi = \varphi_c$, $\theta \in (0, \theta_c)$. A point $g$ on the geodesic $\overline{bc}$ is described by the condition that $a$, $c$, and $g$ lie in the same plane, i.e., $\det(acg) = 0$. This gives the relation

$$\coth\theta_g = A\sin\varphi_g + B\cos\varphi_g = \frac{\sin(\varphi_c - \varphi_g)}{\sin\Phi}\coth\theta_b + \frac{\sin(\varphi_g - \varphi_b)}{\sin\Phi}\coth\theta_c$$

with

$$A = \frac{\coth\theta_c \cos\varphi_b - \coth\theta_b \cos\varphi_c}{\sin\Phi} \quad , \quad B = \frac{\coth\theta_b \sin\varphi_c - \coth\theta_c \sin\varphi_b}{\sin\Phi} \ .$$

The oriented area element on $\mathbf{H}^2$ is given by $\sinh\theta\, d\theta \wedge d\varphi$, and thus the oriented area $\Omega$ of the triangle $\Delta(abc)$ is given by

$$\Omega = \int_{\varphi_b}^{\varphi_c}\left(\int_0^{\theta_g(\varphi)} \sinh\theta\, d\theta\right) d\varphi = C_{bc} + C_{cb} - \Phi\ ,$$

$$C_{jk} = \arcsin\left(\frac{|\sin\Phi|}{\sin\Phi} \cdot \frac{\coth\theta_j - \coth\theta_k \cos\Phi}{\sqrt{\coth^2\theta_j + \coth^2\theta_k - 2\coth\theta_j \coth\theta_k \cos\Phi - \sin^2\Phi}}\right)\ .$$

As before, we introduce $\alpha_{jk} = (\coth\theta_j - \coth\theta_k \cos\Phi) - i\sin\Phi/\sinh\theta_k$ to obtain

$$-i\,e^{iC_{jk}} = \sin C_{jk} - i\cos C_{jk} = \frac{|\sin\Phi|}{\sin\Phi} \cdot \frac{\alpha_{jk}}{|\alpha_{jk}|}\ ,$$

Another elementary computation gives

$$\alpha_{cb} = \frac{e^{-i\varphi_b}(1 - z_b\bar{z}_c)(z_b - z_c)}{2|z_b z_c|} \quad \& \quad \alpha_{bc} = \frac{e^{i\varphi_c}(1 - z_b\bar{z}_c)(\bar{z}_c - \bar{z}_b)}{2|z_b z_c|}\ .$$

It follows that

$$e^{i\Omega} = e^{i(A_{bc}+A_{cb}-\Phi)} = \frac{(1-z_b\bar{z}_c)^2}{|1-z_b\bar{z}_c|^2}\ .$$

Since $|z_b\bar{z}_c| < 1$, it follows that $\mathrm{Arg}(1 - z_b\bar{z}_c) \in (-\pi/2, \pi/2)$. Hence we can extract the square from the argument and we find for the oriented area :

$$\Omega = 2\,\mathrm{Arg}(1 - z_b\bar{z}_c) = \mathrm{Arg}\left(\frac{1 - z_b\bar{z}_c}{1 - \bar{z}_b z_c}\right)\ .$$

As a consequence we find that the area is always smaller than $\pi$, in accordance with the fact that the area is the deficit angle of the triangle, which obviously cannot exceed $\pi$. In order to express this result in a way which does not presuppose that $a$ is the north pole, we note that we have the following identities when $a$ *is* the north pole :

(3) $$\langle a,b\rangle_L = h_b\ ,\quad \langle a,c\rangle_L = h_c\ ,\quad \langle b,c\rangle_L = h_b h_c - \Re(\overline{w}_b w_c)$$
$$\det\nolimits_3(abc) = \Im(\overline{w}_b w_c)\ .$$

Using $z = w/(1+h)$ we find $1 - z_b\bar{z}_c = \dfrac{1 + h_b + h_c + h_b h_c + \overline{w}_b w_c}{(1+h_b)(1+h_c)}$, and thus we obtain the following lemma.

**Lemma.** *The oriented area of the hyperbolic triangle $\Delta(abc)$ is given by*

(4) $$\Omega = 2\,\mathrm{Arg}(1 - z_b\bar{z}_c) = \mathrm{Arg}\left(\frac{1 - z_b\bar{z}_c}{1 - \bar{z}_b z_c}\right) \quad \text{if $a$ is the north pole}$$

(5) $$= 2\,\mathrm{Arg}\Big(1 + \langle a,b\rangle_L + \langle b,c\rangle_L + \langle a,c\rangle_L + i\det\nolimits_3(abc)\Big)\ .$$

*Remark.* Formula (4) can be found in [Pe,p.63] (see also [Be]). Formula (5) can be found in [Ma] and [Ur], who calculated it in the context of relativistic addition of velocities.

We now turn our attention back to the area of a triangle determined by the midpoints of its sides. Without loss of generality we may assume that $\gamma$, the midpoint of the side $\overline{ab}$, is the north pole, which implies in particular that $z_a = -z_b$, $w_a = -w_b$, and $h_a = h_b$. We thus can cut the triangle $\Delta(abc)$ into two triangles $\Delta(\gamma bc)$ and $\Delta(\gamma ca)$ which both satisfy the condition that one corner ($\gamma$) is the north pole. Applying the previous lemma thus gives us for the surface $\Omega$ of the triangle $\Delta(abc)$ the formula :

$$\Omega = 2\operatorname{Arg}(1 - z_b\overline{z}_c) + 2\operatorname{Arg}(1 - z_c(-\overline{z}_b))$$
$$(6) \qquad = 2\operatorname{Arg}\Big((1 - z_b\overline{z}_c)(1 + z_c\overline{z}_b)\Big) \ .$$

This formula gives the area of the triangle in terms of two of its corners ($b$ and $c$) and one midpoint (the north pole $\gamma$). Our next step will be to compute the remaining midpoints in terms of the corners, and then to invert this map. That way we will get the area in terms of the three midpoints.

Putting $\sigma_\beta = w_c + w_a = w_c - w_b$, $\sigma_\alpha = w_c + w_b$, and $\tau = h_a + h_b = h_b + h_c$, we find for the midpoints $\alpha$ and $\beta$ of the sides $\overline{bc}$ and $\overline{ca}$ respectively, the coordinates $w_\epsilon = \sigma_\epsilon/\sqrt{\tau^2 - |\sigma_\alpha|^2}$ and $h_\epsilon = \tau/\sqrt{\tau^2 - |\sigma_\alpha|^2}$, $\epsilon = \alpha, \beta$. A direct computation shows that $\det(\alpha\beta\gamma)^2 = (\Im \overline{w}_\beta w_\gamma)^2 = \dfrac{(\Im w_b \overline{w}_c)^2}{\tau^2 + (\Im w_b \overline{w}_c)^2}$. Hence for the midpoints we always have the constraint $|\det(\alpha\beta\gamma)| < 1$. Since on the hyperbolic plane one can certainly have a triplet $\alpha\beta\gamma$ with a determinant bigger than 1, this is a real constraint and shows that not all triplets $\alpha, \beta, \gamma \in \mathbf{H}^2$ can be obtained as the midpoints of the sides of a triangle.

In order to find the inverse of the map $(z_b, z_c) \mapsto (z_\alpha, z_\beta)$ and to check that the image is completely determined by the condition $|\det(\alpha\beta\gamma)| < 1$, we proceed as follows. We introduce $s_c = \frac{1}{2}(w_\beta/h_\beta + w_\alpha/h_\alpha) = w_c/\tau$ and $s_b = \frac{1}{2}(w_\alpha/h_\alpha - w_\beta/h_\beta) = w_b/\tau$. We thus have $w_j = (h_b + h_c)s_j$ and $h_j^2 - |w_j|^2 = 1$ ($j = b, c$). Solving these equations for $h_b$ and $h_c$ (both$\geq 1$) gives :

$$h_b = \frac{1 + |s_b|^2 - |s_c|^2)}{\sqrt{(1 + |s_b|^2 - |s_c|^2)^2 - 4|s_b|^2}} \quad , \quad \tau = \frac{2}{\sqrt{(1 + |s_b|^2 - |s_c|^2)^2 - 4|s_b|^2}} \ .$$

Since $w_j = \tau s_j$, and using $(1 + |s_b|^2 - |s_c|^2)^2 - 4|s_b|^2 = \dfrac{1 - (\Im \overline{w}_\alpha w_\beta)^2}{h_\alpha^2 h_\beta^2}$, we find the unique solutions

$$h_b = \frac{h_\alpha h_\beta - \Re(w_\beta \overline{w}_\alpha)}{\sqrt{1 - (\Im(w_\beta \overline{w}_\alpha))^2}} \quad , \quad h_c = \frac{h_\alpha h_\beta + \Re(w_\beta \overline{w}_\alpha)}{\sqrt{1 - (\Im(w_\beta \overline{w}_\alpha))^2}}$$
$$w_b = \frac{w_\alpha h_\beta - w_\beta h_\alpha}{\sqrt{1 - (\Im(w_\beta \overline{w}_\alpha))^2}} \quad , \quad w_c = \frac{w_\beta h_\alpha + w_\alpha h_\beta}{\sqrt{1 - (\Im(w_\beta \overline{w}_\alpha))^2}} \ .$$

It follows that the map $(z_b, z_c) \mapsto (z_\alpha, z_\beta)$ is indeed bijective onto $|\det(\alpha\beta\gamma)| < 1$, $\gamma$ being the north pole. Converting this into $z = w/(1+h)$ coordinates yields

$$(1 - z_b \overline{z}_c)(1 + \overline{z}_b z_c) = 2\frac{\tau + i\Im(\overline{w}_b w_c)}{(1 + h_b)(1 + h_c)} \ ,$$

and thus we find from formula (6) for the area :

$$\Omega = 2\operatorname{Arg}\Bigl(\sqrt{1-(\Im(w_\beta \overline{w}_\alpha))^2} + i\Im(w_\beta \overline{w}_\alpha)\Bigr) \ .$$

Since the third midpoint is the north pole, we can apply formula (3) to show that $\Im(w_\beta \overline{w}_\alpha) = \det_3(\alpha\beta\gamma)$. We thus have proven the following proposition.

**Proposition.** *Let $\alpha$, $\beta$, and $\gamma$ be three points on the hyperbolic plane. Then there exists a (unique) triangle $\Delta(abc)$ whose midpoints of the three sides are the given points $\alpha$, $\beta$, and $\gamma$ if and only if $|\det(\alpha\beta\gamma)| < 1$. In that case the oriented area $\Omega \in (-\pi, \pi)$ of $\Delta(abc)$ is given by*

$$\mathrm{e}^{i\Omega/2} = \sqrt{1 - \det(\alpha\beta\gamma)^2} + i\det(\alpha\beta\gamma) \ .$$

## Acknowledgements

I thank S. Berceanu and M. Bordemann for stimulating discussions and helpful remarks.

Laboratoire Paul Painlevé, U.M.R. CNRS 8524 & UFR de Mathématiques; Université de Lille I; F-59655 Villeneuve d'Ascq Cedex; France

*E-mail address*: `Gijs.Tuynman[at]univ-lille1.fr`